 

\baselineskip=14pt
\parskip=10pt

\magnification=\magstephalf

\def\W{{\cal W}}

\def\1{{\overline{1}}}
\def\2{{\overline{2}}}
\parindent=0pt
\overfullrule=0in

\def\frac#1#2{{#1 \over #2}}
\centerline
{\bf
How to Answer Questions of the Type:
}
\centerline
{\bf
If you toss a coin n times, how likely is HH to show up more than HT? 
}
\bigskip
\centerline
{\it Shalosh B. EKHAD and Doron ZEILBERGER}

\qquad\qquad\qquad {\it Dedicated to Dr. Tamar Zeilberger}

{\bf Preface}

On March 16, 2024, Daniel Litt, in an X-post [L] (see also [C]), proposed the following brainteaser.

``Flip a fair coin 100 times. It gives a sequence of heads (H) and tails (T). For each HH in the sequence of flips, Alice gets a point; for each HT, 
Bob does, so e.g. for the sequence THHHT Alice gets 2 points and Bob gets 1 point. Who is most likely to win?"

We show the power of symbolic computation, in particular the (continuous) Almkvist-Zeilberger algorithm, to answer this, and far more general, questions of this kind.
Everything is implemented in our Maple package {\tt Litt.txt} .

{\bf The Maple package Litt.txt}

Our Maple package, {\tt Litt.txt}, can not only answer the original question, but far more general ones, as we will see below. 
First make sure that you have Maple, and that the Maple package, downloadable from

{\tt https://sites.math.rutgers.edu/\~{}zeilberg/tokhniot/Litt.txt} \quad,

resides in your working directory.

Now start a Maple session/worksheet, and load the package, by typing:

{\tt read `Litt.txt`:}

To get the answer whether Alice or Bob are more likely to win the bet, and their respective chances, type:

{\tt WhoWonN(2,$\{$ [1,1]$\}$,$\{$[1,2]$\}$,100);} \quad,

and get right away:

{\tt [false, 0.4576402592, 0.4858327983]} \quad,

meaning that it is {\it false} that Alice is more likely than Bob to win the bet,
and Alice's chance is  $0.4576402592\dots$, while Bob's chance is  $0.4858327983\dots$. It also implies that with probability
$1- 0.4576402592\dots-0.4858327983\dots= 0.0565269425\dots$ it would end with a tie.

How about if you toss $200$ times?, then type

{\tt WhoWonN(2,$\{$ [1,1] $\}$,$\{$[1,2] $\}$,200);} \quad,

and get:

{\tt [false, 0.4700634942, 0.4900044947]} \quad.

Of course, as the number of coin tosses goes to infinity, the chances both converge to $\frac{1}{2}$.
To see why, note that since $11$ and $12$ are both of of same length (namely $2$), the {\bf expected} numbers of occurrences  are the {\bf same} (namely $(n-1)\cdot (\frac{1}{2})^2=(n-1)/4$).

It also happens (see below) that  Bob is {\bf always} more likely to win than Alice for any {\bf finite} $n$.

How about rolling a fair (standard, six-faced) die? If Eve rolls such a die $200$ times, and Alice bets that there are more $11$ than $12$ then, to find out the answer type

{\tt WhoWonN(6,$\{$[1,1]$\}$,$\{$[1,2]$\}$,200);} \quad,

getting {\tt [false, 0.4292455296, 0.4486924385]}, so once again Bob is more likely to win. If it is $11$ versus $23$ then

WhoWonN(6,$\{$ [1,1] $\}$, $\{$ [2,3] $\}$,200);

gives {\tt [false, 0.4346673623, 0.4527404645]}, and once again Bob is more likely to win, but both their chances are better than before, and the
probability of a tie decreased.

But Alice and Bob can also bet on {\it sets} of consecutive strings, as long they are all of the same length. For example

{\tt WhoWonN(6,$\{$[1,1,1],[2,2,2]$\}$,$\{$[1,2,3]$\}$,100);} \quad,

returns {\tt [true, 0.4163070114, 0.1955648145]}, that tells you, not very surprisingly, that if Alice bets that the number of $111$ plus the number of $222$ exceeds the number of $123$, and Bob bets
the opposite, she  wins with probability  $0.4163070114$. On the other hand, entering

{\tt WhoWonN(6,$\{$ [1,1,1],[2,2,2] $\}$,$\{$[1,2,3],[3,2,1] $\}$,100);} \quad, 

gives {\tt [false, 0.3828838919, 0.4121794361]}, and once again Alice is less likely to win
if she predicts that the total number of occurrences of members of $\{111,222\}$ would exceed the total number of occurrences of members of $\{123,321\}$.
Also Alice's and Bob's chances are, respectively,  $0.3828838919$ and $0.4121794361$.

{\bf Using symbolic computation to get numerical answers}

Procedure {\tt WhoWonN(m,A,B,K)} described in the previous section is {\it numeric}, i.e. both its inputs and outputs are {\bf numbers}.
Later on we will also describe more sophisticated procedures that output {\it symbolic} answers, but even for
the simple-minded {\tt WhoWonN}, the most efficient way is via {\it symbol-crunching}.

Of course, a purely brute-force approach won't go very far.
We can generate all  $m^K$ $K$-letter words, and for each of these words count the number of occurrences of consecutive strings that belong to $A$
and those that belong to $B$ and for each find out whether it is an Alice-win, Bob-win, or a tie.

The key to the {\it numeric} procedure {\tt WhoWonN(m,A,B,K)}, as well as the forthcoming {\it symbolic} procedure
{\tt WhoWon(m,A,B,K)},  is to get an {\bf efficient} way to compute the generating function, in the three formal
variables $x$, $a$, $b$, such  that in its Maclaurin expansion the coefficient of $x^n a^i b^j$ would be the number
of $n$-letter words in the alphabet $\{1, \dots, m\}$ with $i$ occurrences, as consecutive subwords, of members of $A$,
and $j$ occurrences, as consecutive subwords, of members of $B$. It is convenient to use the language of {\it weight-enumerators}.

{\bf Computing the Weight-Enumerator with Positive Thinking}

Fix our alphabet to be $\{1, \dots, m\}$, and assume, as we do in this paper, and Maple package, that all the members of $A$ and $B$ are
of the same length, and let that length be $k$. Note that if this is not the case, one can easily transform it to this case by
adding all the possible ways to complete it to sets of $k$-letter words. 

Let $x$, $a$ and $b$ be abstract, `formal' variables, and define the {\bf weight}, $Weight(w)$,  of a word $w=w_1 \dots w_n$, in the alphabet $\{1,2,\dots, m\}$ by
$$
Weight(w):=\, x^n \cdot \prod_{i=1}^{n-k+1} a^{\chi(w_i w_{i+1} \dots w_{i+k-1} \in A)} \cdot b^{\chi(w_i w_{i+1} \dots w_{i+k-1} \in B)}  \quad,
$$
where for any statement $P$, that is either true or false,  $\chi(P)=1$ if $P$ is true, and  $\chi(P)=0$ if $P$ is false. Note that this is the same
as
$$
x^{LengthOfw} \cdot a^{NumberOfOccurrencesOfSubwords \,\, \in \,\, A} \cdot b^{NumberOfOccurrencesOfSubwords \, \, \in \,\, B} \quad.
$$
For example if $A=\{111\}$ and $B=\{222\}$, then
$$
Weight(1111222211)=a^2b^2 \quad .
$$

For any set of words $S$, let $Weight(S)$ be the sum of the weights of the words in $S$. 

We are interested in the {\it weight-enumerator} of the set of {\bf all} words in the alphabet $\{1, ..., m\}$, including the empty word $\phi$.
Let's call this (infinite) set $\W_m$. For $i \geq 0$, let $W_m^{(i)}$ be the set of $i$-letter words.

For each  $v \in W_m^{(k-1)}$, let $\W_m (v)$ be the (infinite) set of words, of length $\geq k-1$, that start with $v$.

For each  $v=v_1 \dots v_{k-1} \in \W_m^{(k-1)}$, we have:
$$
Weight(\W_m(v)) \, = \, Weight(v)+ x \sum_{i=1}^{m} \,  a^{\chi(v_1 \dots v_{k-1} i \in A)} \cdot b^{\chi(v_1 \dots v_{k-1} i \in B)}  \cdot Weight(\W_m(v_2 \dots v_{k-1} i)) \quad .
$$

This gives us a linear system of $m^{k-1}$ equations with $m^{k-1}$ unknowns, that Maple can {\tt solve}. Once we have them, we have
$$
Weight(\W_m)= \sum_{i=0}^{k-2} Weight (\W_m^{(i)})+\sum_{v \in \W_m^{(k-1)}} Weight(\W_m(v)) \quad.
$$
Let's call this grand-generating function $F(x;a,b)$. Note that, thanks to {\it Cramer's law}, this is a {\bf rational function} in the $3$
variables $x,a,b$.

{\bf Computing the Weight-Enumerator with Negative Thinking}

An alternative, more efficient, way of computing $F(x;a,b)$ (with far fewer equations and unknowns) is via the {\bf Goulden-Jackson Cluster method},
that is a {\it negative} approach, using {\it inclusion-exclusion}. See [NZ] for a lucid exposition.

This is implemented in procedure {\tt GFwtE(m,A,B,a,b,x)}. In particular, to get the generating function for the motivating example of the title, type

{\tt GFwtE(2, $\{$ [1,1]  $\}$, $\{$ [1,2]  $\}$,a,b,x);} \quad,

getting
$$
-\frac{a x -x -1}{a \,x^{2}-b \,x^{2}-a x -x +1} \quad.
$$

For old time's sake, we also implemented the slower, positive, approach, and the function call is {\tt GFwtEold(m,A,B,a,b,x);}

If $m^{k-1}$ is large, then of course, the negative approach is {\bf much} faster.

For example

{\tt GFwtE(6,$\{$[1\$6] $\}$, $\{$[1,2\$5]$\}$,a,b,x);}

takes $0.01$ seconds,  yielding that the generating function of words in $\{1,2,3,4,5,6\}$ according to length of word, and number of occurrences of $111111$ and $122222$ is
$$
(-a \,x^{5}-a \,x^{4}+x^{5}-a \,x^{3}+x^{4}-a \,x^{2}+x^{3}-a x +x^{2}+x +1) \cdot
$$
$$
(a b \,x^{10}+a b \,x^{9}-a \,x^{10}-b \,x^{10}+a b \,x^{8}-a \,x^{9}-b \,x^{9}+x^{10}+a b \,x^{7}-a \,x^{8}-b \,x^{8}+x^{9}-a \,x^{7}-b \,x^{7}+x^{8}+5 a \,x^{6}
$$
$$
-b \,x^{6}+x^{7}+5 a \,x^{5}-4 x^{6}+5 a \,x^{4}-5 x^{5}+5 a \,x^{3}-5 x^{4}+5 a \,x^{2}-5 x^{3}-a x -5 x^{2}-5 x +1)^{-1} \quad .
$$

Don't even try to do {\tt GFwtEold(6,$\{$[1\$6] $\}$, $\{$[1,2\$5]$\}$,a,b,x));}, it would take for ever!

\vfill\eject

{\bf Computing the probability of Alice Winning}

Once you have the generating function $F(x;a,b)$ it is very easy to compute the first $200$ (or whatever) terms of the sequence ``probability of Alice winning after $n$ rolls''.
Introduce another formal variable $t$ and look at
$$
F(x;t, t^{-1}) \quad .
$$
Now expand it in a Maclaurin expansion with respect to $x$:
$$
F(x;t,t^{-1})= \sum_{n=0}^{\infty} \, f_n(t) \,x^n \quad .
$$
Here $f_n(t)$ are certain {\it Laurent polynomials} in $t$.
Let's break $f_n(t)$ into three parts:

$\bullet$  the one consisting of the positive powers of $t$, let's call it $f_n^{+}(t)$ 

$\bullet$ the constant term (coefficient of $t^0$), let's call it $f_n^0$

$\bullet$ the one consisting of the negative powers of $t$, let's call it $f_n^{-}(t)$. 

We have:

$$
f_n(t)\,=\, f_n^{-}(t) \,+\, f_n^{0} \,+ \, f_n^{+}(t) \quad.
$$
It follows that, after rolling the $m$-sided fair die $n$ times,

$$
Pr(AliceWon)= \frac{f_n^{+}(1)}{m^n} \quad,
$$
$$
Pr(BobWon)= \frac{f_n^{-}(1)}{m^n} \quad,
$$
and
$$
Pr(ItIsAtie)= \frac{f_n^{0}}{m^n} \quad .
$$

This is the idea behind the numeric {\tt WhoWonN} described above. It can easily go up to a few hundred rolls (or tosses), but what about $30000$ of them? This would be hopeless.
More important, we are {\it mathematicians}, not {\it accountants}, can we get an asymptotic expression, in $n$, for these probabilities? Since everything is asymtotically normal,
and the variance of the number of occurrences of a substring belonging to a fixed set 
is proportional to $\sqrt{n}$, it is not hard to see that Alice's and Bob's probabilities, for each specific scenario, have the asymptotics
$$
\frac{1}{2} - \frac{c_{alice}}{\sqrt{n}} \quad,
$$
$$
\frac{1}{2} - \frac{c_{bob}}{\sqrt{n}} \quad,
$$

for some numbers $c_{alice},c_{bob}$, that depend on $m$, $A$, and $B$, and then of course the probability of a tie is asymptotic to

$$
\frac{c_{alice}+c_{bob}}{\sqrt{n}} \quad .
$$

It would be nice to have accurate estimates for these coefficients, and if possible, exact values.

Speaking of a tie, it is {\it d\'ej\`a vu}! Twelve years ago we [EZ1] used the amazing (continuous) {\bf Almkvist-Zeilberger} algorithm [AZ]
to answer these questions in response to a problem raised by Richard Stanley. Let's recall it briefly.

{\bf How to use the continuous Almkvist-Zeilberger algorithm to investigate the probability of a Tie?}

Easy, the quantity of interest, $f_n^{0}$, is the {\bf constant term}, in $t$, of  the coefficient of $x^n$ in $F(x;t,t^{-1})$. 
So the generating function, in $x$, equals the contour integral below:
$$
\sum_{n=0}^{\infty} f_n^{0} \, x^n \, = \, \frac{1}{2\pi i} \int_{|t|=1} \frac{F(x;t,t^{-1})}{t} \, dt \quad .
$$

Calling this quantity $g(x)$, the continuous Almkvist-Zeilberger algorithm (implemented in \hfill\break
{\tt http://sites.math.rutgers.edu/\~{}zeilberg/tokhniot/EKHAD.txt}, but also
included in the present package, Litt.txt, as procedure {\tt AZc}) outputs a {\it differential equation} with polynomial coefficients satisfied by $g(x)$, and
then the computer easily transforms it to a linear {\bf recurrence} equation for the actual coefficients of $g(x)$, namely, for the integer sequence $\{f_n^{0}\}$. Once this recurrence is found (and it is 
very fast!), one can easily compute the first $30000$ terms, or whatever, and get very accurate estimate for $c_{alice}+c_{bob}$. In fact, using
the method of [EZ2] one can get higher order asymptotics.

It turns out that it is just as easy to use the Almkvist-Zeilberger algorithm to get recurrences for $f_n^{+}(1)$ and $f_n^{-}(1)$ enabling the
fast computation of the probabilities of Alice and Bob winning the bet.

To wit

$$
\sum_{n=0}^{\infty} f_n^{+}(1) x^n \, = \, \frac{1}{2\pi i} \int_{|t|=1} \frac{F(x;t,t^{-1})}{1-t} \, dt \quad ,
$$
with a similar expression for the generating function for $f_n^{-}(1)$.

This is implemented in procedure {\tt RECaz(m, A, B, n, N)}, that inputs $m$ (the size of the alphabet), and the sets $A$ and $B$ belonging to Alice and Bob, respectively
and outputs a linear recurrence equation with polynomial coefficients for  $f_n^{+}(1)$. For example

{\tt RECaz(2,$\{$ [1,1] $\}$, $\{$[1,2] $\}$,n,N);}

outputs a certain $7^{th}$-order linear recurrence that enables very fast computation of many terms.

This is all combined in procedure 

{\tt WhoWon(m,A,B,K)} \quad,

that uses $K$ terms of the sequence to estimate the numbers $c_{alice}, c_{bob}$ mentioned above, that tells you that Alice's and Bob's chances of winning after $n$ rolls, are
asymptotic to $\frac{1}{2} - \frac{c_{alice}}{\sqrt{n}}$, and $\frac{1}{2} - \frac{c_{bob}}{\sqrt{n}}$, respectively.

{\bf Some Sample Output}

There are numerous output files in the web-page of this paper:

{\tt https://sites.math.rutgers.edu/\~{}zeilberg/mamarim/mamarimhtml/litt.html} \quad .

Let's just give a few highlights.

For the original problem, type:

{\tt WhoWon(2,$\{$[1,1]$\}$,$\{$[1,2] $\}$,30000)};

getting

[false, 0.423144, 0.141049] \quad,

meaning that it is false that Alice is more likely than Bob to win (i.e. Bob is more likely to win) and that
$$
Pr(AliceWon) \,\,  \asymp \,\, \frac{1}{2}-\frac{ 0.423144}{\sqrt{n}} \quad,
$$
$$
Pr(BobWon) \,\, \asymp \,\, \frac{1}{2}-\frac{ 0.141049 }{\sqrt{n}} \quad .
$$
It turns out (see later) that these are in fact  {\bf exactly}  (see below)  
$\frac{1}{2}-\frac{ 3}{4 \sqrt{\pi} \sqrt{n}}$, and $\frac{1}{2}-\frac{ 1}{4 \sqrt{\pi} \sqrt{n}}$, respectively.

The advantage of Bob over Alice is asymptotic to  $\frac{c_{alice}-c_{bob}}{\sqrt{n}}$. We have, using $K=20000$:

$\bullet$ The best counter bets against Alice's $111$  are Bob's $112$ and (equivalently) $211$, and Bob's advantage is (estimated to be) $\frac{0.598456}{\sqrt{n}}$ \quad .

$\bullet$ The second best counter bets against Alice's $111$  are Bob's $122$ and (equivalently) $221$, and Bob's advantage is (estimated to be) $\frac{0.4886160}{\sqrt{n}}$ \quad .

$\bullet$ The third best counter bet against Alice's $111$  is Bob's $121$  and Bob's advantage is (estimated to be) $\frac{ 0.32572}{\sqrt{n}}$ \quad .

$\bullet$ The fourth best counter bet against Alice's $111$  is Bob's $212$  and  Bob's advantage is (estimated to be) $\frac{ 0.28214}{\sqrt{n}}$ \quad .

$\bullet$ The fifth (and worst) counter bet against Alice's $111$  is Bob's $222$  and  Bob's advantage is $0$, of course, by symmetry.

To see the (complicated) linear recurrences that lead to these estimates, look at the output file:

{\tt https://sites.math.rutgers.edu/\~{}zeilberg/tokhniot/oLitt23.txt} \quad .

{\bf More precise asymptotics}

The recurrences one gets from the Almkvist-Zeilberger algorithm are not, usually, minimal, and hence are not good for deriving asymptotics using the method of [EZ2].
But one can easily conjecture minimal recurrences, and then prove them rigorously, i.e. that the sequences are equivalent, using the {\it Euclidean algorithm}
in the non-commutative algebra of linear recurrence operators with polynomial coefficients. Since it is always possible, we didn't bother.

It turns out that if write:
$$
Pr(AliceWon)(n)= \frac{1}{2} - A_n \quad,
$$
then $A_n$ satisfies a nice recurrence that leads to more precise asymptotics. In particular, for the origianl $HH$ vs $HT$ problem, We have the following two theorems.

{\bf Theorem 1}: In the original Litt game, when you toss a fair coin $n$ times, the probability that $HT$ would beat $TT$ is
$$
\frac{1}{2}- A_n \quad,
$$
where $A_n$ is a solution of the recurrence
$$
-\frac{\left(n +1\right) A_{n}}{8 \left(n +4\right)}+\frac{\left(4 n +7\right) A_{n +1}}{8 n +32}-\frac{\left(5 n +12\right) A_{n +2}}{8 \left(n +4\right)}+\frac{\left(3 n +8\right) A_{n +3}}{4 n +16}-\frac{\left(3 n +10\right) A_{n +4}}{2 \left(n +4\right)}+A_{n +5}
=0 \quad,
$$
subject to the initial conditions
$$
A_{1} = {\frac{1}{2}}, A_{2} = {\frac{1}{4}}, A_{3} = {\frac{1}{8}}, A_{4} = {\frac{1}{8}}, A_{5} = {\frac{3}{32}} \quad.
$$
The asymptotic of $A_n$ starts as
$$
\frac{1}{4 \sqrt{\pi}}
\left (
\sqrt{\frac{1}{n}}+\frac{7 \sqrt{\frac{1}{n}}}{16 n}+\frac{265 \sqrt{\frac{1}{n}}}{512 n^{2}}+\frac{13165 \sqrt{\frac{1}{n}}}{8192 n^{3}}+\frac{3996699 \sqrt{\frac{1}{n}}}{524288 n^{4}}+\frac{377801193 \sqrt{\frac{1}{n}}}{8388608 n^{5}}
+O(n^{-13/2})
\right ) \quad .
$$

{\bf Theorem 2}: In the original Litt game, when you toss a fair coin $n$ times, the probability that $HH$ would beat $HT$ is
$$
\frac{1}{2}- B_n \quad,
$$
where $B_n$ is a solution of the recurrence
$$
\frac{\left(n +1\right) B_{n}}{8 n +40}+\frac{3 B_{n +1}}{8 \left(n +5\right)}-\frac{3 \left(n +4\right) B_{n +2}}{8 \left(n +5\right)}-\frac{\left(n +1\right) B_{n +3}}{4 \left(n +5\right)}-\frac{\left(n +6\right) B_{n +4}}{2 \left(n +5\right)}+B_{n +5}
=0 \quad ,
$$
subject to the initial conditions
$$
B_{1} = {\frac{1}{2}}, B_{2} = {\frac{1}{4}}, B_{3} = {\frac{1}
{4}}, B_{4} = {\frac{1}{4}}, B_{5} = {\frac{3}{16}} \quad .
$$
The asymptotic of $B_n$ starts as
$$
\frac{3}{4 \sqrt{\pi}}
\left (
\sqrt{\frac{1}{n}}+\frac{5 \sqrt{\frac{1}{n}}}{48 n}+\frac{169 \sqrt{\frac{1}{n}}}{512 n^{2}}+\frac{26615 \sqrt{\frac{1}{n}}}{24576 n^{3}}+\frac{2583259 \sqrt{\frac{1}{n}}}{524288 n^{4}}+\frac{242384345 \sqrt{\frac{1}{n}}}{8388608 n^{5}}
+O(n^{-13/2)}
\right ) \quad .
$$

For similar theorems, for more complicated bets, consult the web-page of this paper. 

{\bf Conclusion}: We illustrated the power of symbolic computation to help you play, for profit, Litt-style games.

{\bf Acknowledgment}: Many thanks to Dr. George Spahn (DZ's academic son), and Dr. Tamar Zeilberger (DZ's biological daughter) for telling us about Daniel Litt's X-post [L].

{\bf References}

[AZ] Gert Almkvist and Doron Zeilberger, {\it  The method of differentiating under the  integral sign},  J. Symbolic Computation {\bf 10} (1990), 571-591. \hfill\break
{\tt  https://sites.math.rutgers.edu/\~{}zeilberg/mamarim/mamarimPDF/duis.pdf} \quad .

[C] Roman Cheplyaka, {\it Alice and Bob flipping coins puzzle}, March 18, 2024. \hfill\break
{\tt https://ro-che.info/articles/2024-03-18-alice-bob-coin-flipping}, 

[EZ1] Shalosh B. Ekhad and Doron Zeilberger, {\it Automatic Solution of Richard Stanley's Amer. Math. Monthly Problem \#11610 and ANY Problem of That Type},
The personal Journal of Shalosh B. Ekhad and Doron Zeilberger, Dec. 28, 2012. \hfill\break
{\tt https://sites.math.rutgers.edu/\~{}zeilberg/mamarim/mamarimhtml/rps.html} \quad .\hfill\break
arxiv: {\tt https://arxiv.org/abs/1112.6207} \quad.

[EZ2] Shalosh B. Ekhad and Doron Zeilberger, 
{\it AsyRec: A Maple package for Computing the Asymptotics of Solutions of Linear Recurrence Equations with Polynomial Coefficients}, 
The personal Journal of Shalosh B. Ekhad and Doron Zeilberger, April 4, 2008. \hfill\break
{\tt https://sites.math.rutgers.edu/\~{}zeilberg/mamarim/mamarimhtml/asy.html}

[L] Daniel Litt, {\it X-post}, March 16, 2024. {\tt https://x.com/littmath/status/1769044719034647001}

[NZ] John Noonan and Doron Zeilberger,  {\it The Goulden-Jackson cluster method: extensions, applications, and implementations}, J. Difference Eq. Appl. {\bf 5}(1999), 355-377. \hfill\break
{\tt https://sites.math.rutgers.edu/\~{}zeilberg/mamarim/mamarimhtml/gj.html}

\bigskip
\hrule
\bigskip
Shalosh B. Ekhad, c/o D. Zeilberger, Department of Mathematics, Rutgers University (New Brunswick), Hill Center-Busch Campus, 110 Frelinghuysen
Rd., Piscataway, NJ 08854-8019, USA. \hfill\break
Email: {\tt ShaloshBEkhad at gmail dot com}   \quad .
\bigskip
Doron Zeilberger, Department of Mathematics, Rutgers University (New Brunswick), Hill Center-Busch Campus, 110 Frelinghuysen
Rd., Piscataway, NJ 08854-8019, USA. \hfill\break
Email: {\tt DoronZeil at gmail  dot com}   \quad .
\bigskip
{\bf Exclusively published in the Personal Journal of Shalosh B. Ekhad and Doron Zeilberger, and arxiv.org}
\bigskip
{\bf May 20, 2024} \quad .

\end